\documentclass[11pt]{article}
\usepackage{amsmath}
\usepackage{amsthm}
\advance \topmargin by -\headheight
\advance \topmargin by -\headsep
\evensidemargin \oddsidemargin
\theoremstyle{plain}
\newtheorem{thm}{Theorem}[section]
\newtheorem{prop}[thm]{Proposition}
\newtheorem{lem}[thm]{Lemma}
\newtheorem{cor}[thm]{Corollary}
\newtheorem{conj}[thm]{Conjecture}

\theoremstyle{remark}
\newtheorem{rmk}[thm]{Remark}

\theoremstyle{definition}

\title{\emph{Homology Lens Spaces \\ in Topological $4$-Manifolds}}
\author{Allan L. Edmonds}
\date{}

\begin{document}
\maketitle
\abstract{For a closed $4$-manifold $X^4$ and closed $3$-manifold
$M^3$ we investigate the smallest integer $n$ (perhaps
$n=\infty$) such that $M^3$ embeds in $\#_nX^4$, the connected sum of
$n$ copies of $X^4$.  It is proven that any lens space (or homology lens space)
embeds topologically locally flatly in $\#_2({\mathbf C}P^2\#\ \overline
{{\mathbf C}P}^2)$, in $\#_4 S^2\times S^2$ and in $\#_8
\mathbf{C}P^2$}.

\section{Introduction}

For any closed $4$-manifold $X^4$ and closed $3$-manifold $M^3$ one
can define a simple numerical invariant, the $X^4$-genus of $M^3$, denoted
$g_{X^4}(M^3)$, by saying $g_{X^4}(M^3)\le n$ provided $M^3$ embeds in
$\#_nX^4$, the connected sum of $n$ copies of $X^4$.  (All embeddings are
understood to be locally flat.) It follows that $0\le g_{X^4}(M^3)\le
\infty$ and that $g_{X^4}(M^3)=0$ if $M^3$ embeds in $S^4$, while
$g_{X^4}(M^3)=\infty$ if $M^3$ embeds in no $\#_nX^4$ for any
integer $n$.  Here we understand our manifolds to be topological
manifolds and our embeddings to be locally flat.  There is also an
analogous invariant $g_{X^4}^{DIFF}(M^3)$ where one requires $X^4$
and the embedding to be smooth. (We endow the $3$-manifold with its
unique smooth structure.)  Certainly $g_{X^4}(M^3) \le
g_{X^4}^{DIFF}(M^3)$ for any smooth
$4$-manifold.

For example consider the case $X=S^2\times S^2$.  It is known that
any closed orientable $3$-manifold $M^3$ embeds smoothly in some
$\#_n S^2\times S^2$, so that $g_{S^2\times
S^2}^{DIFF}(M^3)<\infty$.  It will be shown here that
$g_{S^2\times S^2}(L(p,q))\le 4$ for every lens space $L(p,q)$.  We
doubt that the DIFF $S^2\times S^2$-genus of lens spaces is bounded.

It is also known that any lens space $L(p,q)$ embeds smoothly in
some $\#_n \mathbf{C}P^2$, so that
$g_{\mathbf{C}P^2}^{DIFF}(L(p,q))<\infty$.  But not every
$3$-manifold embeds smoothly in some $\#_n
\mathbf{C}P^2$, by gauge-theoretic considerations.  For example, the
Poincar\'e homology $3$-sphere does not embed smoothly in $\#_n
\mathbf{C}P^2$ for any positive integer $n$.  It will be shown here
that
$g_{\mathbf{C}P^2}(L(p,q))\le 8$.  It is not known whether this is a
sharp bound, but it is known  \cite{EL} that 
$g_{\mathbf{C}P^2}(L(8k,1))=  5$ for $k> 2$.  Again, we doubt very
much that
$g_{\mathbf{C}P^2}^{DIFF}(L(p,q))$ is bounded.

One interest in the embedding questions considered here stems from
work of F. Fang \cite{F} who showed that if a
$3$-manifold $M^3$ embeds in $\#_n{\mathbf C}P^2$, then the open
$4$-manifold $M^3\times\mathbf{R}$ admits uncountably many smooth
structures.

Now when a $3$-manifold $M^3$ embeds in a closed, simply connected
$4$-manifold $X^4$, $M^3$ bounds in $X^4$, splitting
$X^4$ into two compact submanifolds $U$ and $V$ with $\partial
U=M=\partial V$.  Taking into account orientations, assuming that
$U$ and $V$ inherit orientations from one on $X^4$ and that $M^3$
(which must then admit an orientation) is oriented, then we can
assume that  $\partial U=M^3$, while $\partial V=-M^3$.  Thus, to
show that
$M^3$ embeds in any particular $X^4$ it suffices to show that $M^3$
and $-M^3$ bound appropriate (preferably simply connected) manifolds
$U$ and $V$ such that $U\cup_MV\simeq X^4$.  To recognize $U\cup_MV$
as $X^4$ it helps to be in the topological category, where one can
apply Freedman's classification in terms of the intersection pairing
(and the Kirby-Siebenmann triangulation obstruction).

Our starting point will be the following result.

\begin{thm} Any 3-dimensional (homology) lens space $L(p,q)$  bounds
a compact, simply connected, topological 4-manifold with
$b_2\le 2$. 
\end{thm}

By exercising due care we can show that a homology lens space always
bounds suitable $4$-manifolds to show that it embeds in certain
relatively small connected sums.

\begin{thm} Any 3-dimensional (homology) lens space $L(p,q)$ embeds
topologically locally flatly in $\#_2({\mathbf C}P^2\#\ \overline
{{\mathbf C}P}^2)$, in 
$\#_4S^2\times S^2$, and in $\#_8{\mathbf C}P^2$.
\end{thm}

If
$p$ is odd or if
$q\equiv
\pm 1
\mod p$, then
$L(p,q)$ actually embeds in
$\#_5{\mathbf {C}}P^2$.  But, as noted above, it is known that
$L(8k,1)$ does not embed in  $\#_4{\mathbf {C}}P^2$ when $k>2$. (See
\cite{EL}.)

\begin{conj} For any simply connected $4$-manifold $X^4$ and for any
(homology) lens space $L(p,q)$, the $X^4$-genus
$g_{X^4}(L(p,q))<\infty$.
\end{conj}

For $X^4$ indefinite, this follows easily from the present work, so
it essentially reduces to considering $X^4$ with a definite
intersection pairing that does not split nontrivially as an
orthogonal sum.  The same conjecture may be posed for any rational
homology 3-sphere in place of the lens space $L(p,q)$.  F. Fang, in
unpublished work, however, has shown that there exist 3-manifolds with
large first betti number $b_1$ that do not embed in any positive definite,
simply connected, 4-manifold.

All the results about topological embeddings of lens spaces derived
in this paper apply equally well to any homology lens space, that is
a
$3$-manifold $M^3$ with $H_1(M^3)$ finite cyclic.   Thus the phrase
``homology lens space $L(p,q)$'' refers to any 3-manifold with
$H_1=\mathbf{Z}_p$ and linking form equivalent to $({(p-q)}/{q})$, or, after
changing orientation, $(q/p)$.

\section{Intersection Pairings and Linking Forms}
\label{intersect}
The
crucial invariants we must study are the linking form associated
with a $3$-manifold, the intersection pairing associated to a
$4$-manifold bounded by the $3$-manifold, and the relationship
between the two.
\subsection{Intersection numbers and linking numbers} We refer to the
classic book of Seifert and Threlfall [1934], Sections 73-77, for
generalities and basic definitions of intersection numbers and
linking numbers, in the context of polyhedral manifolds.

\subsection{Intersection pairings and linking forms} By an {\it
abstract intersection pairing} we understand a finitely generated
free abelian group $F$ together with a symmetric bilinear mapping
${\mathcal S}:F\times F\to {\mathbf Z}$ (${\mathcal S}$ for
``Schnittzahlen'').  We will only be concerned with nondegenerate
intersection pairings, such that the associated adjoint homomorphism
${\rm ad}{\mathcal S}:F\to {\rm Hom}(F,{\mathbf Z})$ has nonzero
determinant, or equivalently has finite index image.  It is often
convenient to describe such pairings by square integer matrices that
give the adjoint ${\rm ad}{\mathcal S}$ with respect some basis of
$F$ and the corresponding dual basis for ${\rm Hom}(F,{\mathbf
Z})$.  Our main geometric example of such an intersection pairing is
the second homology of a simply connected 4-manifold with connected
boundary, where the boundary is a rational homology 3-sphere, under
the usual intersection number.  On the algebraic side any symmetric
integer matrix with nonzero determinant determines such an
intersection pairing.  And all such algebraic intersection pairings
can be realized by smooth 4-manifolds, by attaching 2-handles to the
4-ball along a suitable framed link.

By an {\it abstract linking form} we understand a finite abelian
group $G$ together with a symmetric bilinear mapping ${\mathcal
V}:G\times G\to {\mathbf Q}/{\mathbf Z}$ ({$\mathcal V$} for
``Verschlingungszahlen'').  We will only be concerned with
nonsingular linking forms, such that the associated adjoint
homomorphism ${\rm ad}{\mathcal V}:G\to {\rm Hom}(G,{\mathbf
Q}/{\mathbf Z})$ is an isomorphism.  One can describe such a linking
form by a suitable matrix of rational numbers, by choosing, say,
generators for $G$, and indicating the pairings of pairs of these
generators.  See also the end of this section.

Our main geometric example of such a linking form is the linking
form of a 3-manifold.  Let $M^3$ be a closed oriented 3-manifold and
let
$T_1(M^3)$ denote the torsion subgroup of the first homology
$H_1(M^3)$, with integer coefficients. We define the classical
linking form ${\mathcal V}_{M^3} :T_1(M^3)\times T_1(M^3)\to
{\mathbf Q}/{\mathbf Z}$  as follows:  Suppose that
$\alpha, \beta \in T_1(M^3)$.  Represent $\alpha$ and
$\beta$ by disjoint 1-cycles, or even simple closed curves
$A$ and $B$.  There is a positive integer
$n$ such that $n\beta=0$ in $H_1(M^3)$.  Thus there is a 2-cycle $C$
with $\partial C = nB$.  Then
${\mathcal V}_{M^3}(\alpha , \beta) = A\cdot C/n$ in ${\mathbf
Q}/{\mathbf Z}$. One argues that the linking form is well-defined,
independent of all the choices made in the course of its
definition.  The linking form is symmetric and is  nonsingular, in
the sense that the associated homomorphism $T_1(M^3)\to {\rm Hom}
(T_1(M^3),{\mathbf Q}/{\mathbf Z})$ is an isomorphism, by Poincar\'e
Duality and Universal Coefficients.

\subsection{Presentation of linking forms}

Suppose that $M^3=\partial W^4$, where $W^4$ is a compact oriented
{\it simply connected} 4-manifold.  Every closed oriented 3-manifold
bounds such a 4-manifold.  Consider the homology long exact sequence
of the pair
$(W^4,M^3)$.
$$0\longrightarrow H_2(M^3)\buildrel i \over
\longrightarrow H_2(W^4)\buildrel j \over \longrightarrow
H_2(W^4,M^3)\buildrel
\partial \over \longrightarrow H_1(M^3)\longrightarrow 0$$

\begin{lem}\label{pres} If $j(\xi_i)=n_i\alpha_i$ (i=1,2)
 where $n_1n_2\ne 0$, then ${\mathcal V}_{M^3} (\partial\alpha_1,
\partial\alpha_2)=  \frac{-1}{ n_1n_2}{\mathcal S}(\xi_1,\xi_2)$ in
${\mathbf Q}/{\mathbf Z}$, where $\mathcal S$ denotes the
intersection pairing on $H_2(W^4)$.  Also, if $u\in H_2(M^3)$ and
$\eta\in H^2(W^4)=H_2(W^4,M^3)$, then $\partial(\eta)\cdot
u=\eta(i(u))$.
\end{lem}

The last condition is null in the case the boundary is a rational
homology sphere or, equivalently, the intersection pairing in
nondegenerate. One says that
$H_2(W^4)$ together with the intersection pairing
$\mathcal S$ {\it presents}
$H_1(M^3)$. In general a nondegenerate intersection pairing
$(F,{\mathcal S})$ {\it presents} a linking form $(G,{\mathcal V})$
if there is a short exact sequence $$ 0\longrightarrow F\buildrel
{\rm ad}{\mathcal S} \over \longrightarrow {\rm Hom}(F,{Z})\buildrel
\partial
\over
\longrightarrow G\buildrel
 \over \longrightarrow 0$$  such that if ${\rm ad}{\mathcal
S}(\xi_i)=n_i\alpha_i$ (i=1,2)
 where $n_1n_2\ne 0$, then ${\mathcal V} (\partial\alpha_1,
\partial\alpha_2)=  \frac{-1}{n_1n_2}{\mathcal S}(\xi_1,\xi_2)$ in
${\mathbf Q}/{Z}$.  Another way of writing this (compare Turaev
[1984], Section 3) is 
${\mathcal V}(\partial\alpha_1, \partial\alpha_2)=  -{\mathcal
S}^{-1}(\alpha_1,\alpha_2)$ in
${\mathbf Q}/{\mathbf Z}$.  Here ${\mathcal S}^{-1}$ denotes the
restriction of the rational bilinear pairing whose associated
homomorphism is $(ad{\mathcal S}\otimes {\mathbf Q})^{-1}:{\rm
Hom}(F,{\mathbf Q})={\rm Hom}(F,{\bf Z})\otimes{\mathbf Q}\to
F\otimes {\mathbf Q}$. 

One consequence of Wall's analysis \cite{W} of
linking forms is that any abstract linking form is presented by some abstract
intersection pairing.  For example, $(q/p)$ is presented by a matrix whose size
depends on a certain  continued fraction decomposition of
$q/p$.  Geometrically this corresponds to expressing a lens space as
the boundary of an appropriate plumbing manifold, whose second betti
number depends on the length of the corresponding continued
fraction.  One of our goals is to bound the rank of such
presentations.

 
\subsection{Topological realization of algebraic presentations}
 Here we
describe known conditions for translating a presentation into a
${4}$-manifold with given boundary. 


Given an abstract intersection pairing, described, say, by a
symmetric integer matrix with nonvanishing determinant, then this
intersection pairing can easily be realized by a compact simply
connected smooth 4-manifold obtained by attaching 2-handles to the
4-ball along any framed link whose linking matrix is the given
symmetric integer matrix.   In this way one hardly controls the
corresponding boundary, except to say that its linking pairing is
determined as above.  It turns out that in the topological category
one can say much more, actually prescribing the boundary in advance.
S. Boyer \cite{B} and Stong \cite{St} have independently proven the following
result, which extends Freedman's original realization result for closed
simply connected 4-manifolds.

\begin{thm}\label{real}  If the geometric linking form
$(H_1(M^3),{\mathcal V})$ is presented by an abstract intersection
pairing
$(F,{\mathcal S})$, then
$M^3$ is the boundary of a simply connected topological 4-manifold
$X^4$ with $H_2(X^4)=F$ and ${\mathcal S}$ as intersection pairing.
\end{thm}

One corollary of this result is that if a given 3-manifold
$N^3$ bounds a simply connected 4-manifold $Y^4$, then any other
3-manifold $M^3$ with the same linking form bounds a topological
4-manifold
$X^4$ with the same intersection pairing as $Y^4$.  Going further,
the work of Boyer actually characterized when two simply connected
4-manifolds with given boundary are homeomorphic.

So, in the topological category, the question of what kinds of
simply connected 4-manifolds  have a given boundary is in fact
reduced to a purely algebraic one about existence of suitable
presentations of linking forms.

\begin{rmk}\label{ks} When the intersection pairing on $F$ has odd
type, the work of Boyer and of Stong shows that both a 0 and a
nonzero Kirby-Siebenman stable triangulation obstruction in
$H^4(X^4,M^3;{\mathbf Z}_2)$ can be realized.
\end{rmk}

Both Boyer and Stong dealt with general closed oriented 3-manifolds,
not just rational homology spheres, as considered here.  In the
present case a rather simpler proof is available, which we sketch
for the reader's convenience (cf. Boyer\cite{B}, Section 8).

\begin{proof}[Sketch  proof of Theorem \ref{real}] Given a
presentation $(F,{\mathcal S})$ of the linking form on
$H_1(M^3)$, we can realize $(F,{\mathcal S})$ by a framed link in
the 3-sphere.  Attaching 2-handles to the 4-ball along this framed
link produces a smooth, compact, simply connected 4-manifold $V^4$
with intersection pairing 
$(F,{\mathcal S})$.  The boundary $\partial V^4 = N^3$ is a
3-manifold with a linking form equivalent to that of the given
3-manifold $M^3$.  Passing to the dual handle decomposition, we see
that $V^4$ can be described as being obtained from $N^3\times I$ by
attaching 2-handles along a framed link in $N^3=N^3\times \{ 1\}$,
and then capping off with a 4-handle.  We can choose a framed link
in $M^3$ that mirrors this link in $N^3$, in the sense that the
elements of the link represent corresponding elements in first
homology and all linking numbers and framings agree with those in
$N^3$.  If we add 2-handles to $M^3\times I$ along this framed link
in $M^3\times \{ 0\}$, we obtain a compact, smooth 4-manifold $W^4$,
with one boundary component $M^3$ and the other boundary component a
homology 3-sphere
$\Sigma^3$.  By Freedman we can cap off $\Sigma^3$ with a compact
contractible topological 4-manifold
$\Delta^4$.  In particular $X^4=W^4\cup\Delta^4$ is a compact simply
connected, topological 4-manifold with boundary $M^3$ and
intersection pairing equivalent to $(F, {\mathcal S})$.
\end{proof}

\section{Minimal presentations of $(q/p)$}

In this section we will derive smallest possible presentations of the
indecomposable linking forms $(q/p)$.  

\subsection{Rank $1$ forms} Here we determine which rank $1$ linking
forms are presented by  rank $1$ intersection pairings. 

A rank 1 linking form on ${\mathbf Z}/p$ can be described by a $1
\times 1$ matrix $(q/p)$, where $q$ is prime to $p$ and
$q$ is well-defined up adding a multiple of $p$ and multiplying by a
square of a unit mod $p$.  It is the linking form of the lens space
$L(p,q)$.  Such a linking form can always be presented by some
intersection pairing.  One way to do this is to develop a continued
fraction expansion of $p/(p-q)$, as in Hirzebruch et al \cite{HNK}. Say it is
$[a_1,\dots , a_n]$.  This defines a plumbing 4-manifold
$P^4[a_1,\dots , a_n]$, which has an intersection pairing of rank
$n$ and has oriented boundary $L(p,q)$.  As $p$ and
$q$ vary the rank $n$ does not stay bounded.  We seek a way of
controlling the rank.

The
$1\times 1$ intersection pairing
$(p)$ is positive definite and realizes the linking form
$(-1/p)$.  We state this as follows:

\begin{thm}  The linking form $(q/p)$ is realized by a rank 1 matrix
$\left[ a\right]$ if and only if $\mp q$ is a quadratic residue mod
$p$ (and $a=\pm p$).
\end{thm}


\subsection{Rank $2$ forms.} We begin the study of rank 2
presentations with a general realization statement, which provided
the original starting point for this paper.

\begin{thm} If $p$ and $q$ are relatively prime integers, then the
abstract linking form
$(q/p):{\mathbf Z}/p\times {\mathbf Z}/p\ \to {\mathbf Q}/{\mathbf
Z}$ is presented by a non-degenerate rank $2$ abstract intersection
pairing
${\mathcal S}:{\mathbf Z}^{2}\times {\mathbf Z}^{2} \to {\mathbf Z}$
of odd type.
\end{thm}

\begin{proof} For the purposes of the proof we may assume that
$p$ is positive.  Note that we may also replace $q$ by $-q$ if we
wish. If $b$ and $d$ are integers, then the matrix
$$ G=\left[
\begin{matrix} q/p & b\\ b & d
\end{matrix}
\right]
$$ understood mod ${\bf Z}$ also gives the linking form
$(q/p)$.  If one can choose the integers $b$ and $d$ so that $\det
G=\pm1/p$, then $G^{-1}$ is integral and ${\mathcal S}=-G^{-1}$
presents $(q/p)$.  

Now $\det G=dq/p -b^2$, so this amounts to solving the equation
$dq-pb^2=\pm 1$ for $b$ and $d$.  Actually we have a little more
freedom, since $q$ is only defined modulo $p$.

One of $q$, $-q$, $p+q$, or $3p+q$  must be congruent to 3 mod 4. 
Thus we may assume that
$q\equiv 3$ mod 4.  Now consider the arithmetic progression
$q+4np$, $n=1, 2, \ldots$.  By Dirichlet's theorem on primes in an
arithmetic progression (see \cite{IR}, for example), some
$q+4np$ is prime.  Therefore we can assume that
$q$ is a prime congruent to 3 mod 4.  But then, since $-1$ is not a
square mod $q$, either $p$ or $-p$ must be a square mod $q$ and the
proposition follows.
\end{proof}

\begin{rmk} Note that we have shown that $(\pm q/p)$ is presented by
the intersection pairing
$$ {\mathcal S}=-G^{-1}=
\left[
\begin{matrix} -dp & bp \\ bp & -q
\end{matrix}
\right]
$$ and in particular has a diagonal entry that is negative and odd.
\end{rmk}

\section{First Topological Applications}

Here we combine the theorem of Boyer and Stong with the algebraic
result of the preceding section to find small coboundaries for lens
spaces and other 3-manifolds.

\begin{cor}\label{b2=1} If $p$ is a positive integer and $q$ is an
integer prime to $p$, then any (homology) lens space $L(p,q)$ is the
boundary of a simply connected topological 4-manifold, with $b_2= 1$ if and
only if
$\pm q$ is a quadratic residue$\mod p$.
\end{cor} Corollary \ref{b2=1} is originally due to O. Saeki
\cite{S}.  Earlier R. Fintushel and R. Stern \cite{FS} studied the
problem of when a lens space bounds such a $4$-manifold
\emph{smoothly} and found both further obstructions and some explicit
constructions.

\begin{cor} If $p$ is a positive integer and $q$ is an integer prime
to $p$, then any (homology) lens space  $L(p,q)$ is the boundary of a simply
connected topological 4-manifold, with $b_2\le 2$.
\end{cor}


We note that the rank 2 presentation matrix above is necessarily of
odd type, as the diagonal entry $\pm q$ is odd.

\begin{cor} If $p$ is a positive integer and $q$ is an integer prime
to $p$, then  any (homology) lens space  $L(p,q)$
 admits a topological embedding in $\#_2({\mathbf{C}}P^2\#\
\overline {\mathbf{C}P}^2)$.
\end{cor}
\begin{proof} We have seen that $L(p,q)$ bounds a simply connected
$4$-manifold of odd type and with $b_2=2$.  The double of such a
$4$-manifold is precisely
$\#_2({\mathbf{C}}P^2\#\ \overline {\mathbf{C}P}^2)$.  This follows
from Freedman's classification of simply connected topological
$4$-manifolds together with the observation that the mod $2$
Kirby-Siebenmann invariant of the double vanishes by additivity.
\end{proof}

\section{Even Intersection Pairings and Applications}

Here we investigate presenting a linking form by an {\em even} intersection
pairing.  We use an algebraic analog of the geometric notions of
``blow-up'' and ``blow-down''.  If
$V$ is an intersection pairing, then we will refer to
$V\oplus \left<+1\right>$ and $V\oplus \left<-1\right>$ as being
blow-ups of
$V$.  If $v\in V$ and $v\cdot v =
\pm 1$, then
$v^\perp =\{ u\in V: u\cdot v = 0 \}$ is an orthogonal summand of
$V$ and we say that $v^\perp$ is obtained from
$V$ by blowing down $v$.  Notice in particular that $V$ and its
blow-ups and blow-downs all present the same linking form, as would
the orthogonal sum of $V$ with any unimodular pairing.

Recall that if $V$ is an intersection pairing, then an element $v\in
V$ is said to be {\it characteristic} if one has
$v\cdot w \equiv w\cdot w \mod 2$ for all $w\in V$.  An element that
is not characteristic is called {\it ordinary}. It is easy to see
that characteristic elements always exist.  (Geometrically, the mod
$2$ choices for characteristic elements correspond to spin
structures on the boundary manifold.)  A key point is that if $v\in
V$ is characteristic, then the induced pairing on $v^\perp$ is of
even type.  Similarly, if $V$ has odd type and
$v\in V$ is not characteristic, then the induced pairing on
$v^\perp$ is of again of odd type.

\begin{prop} The rank 1 linking form $(q/p)$ is presented by a rank
4 intersection pairing of even type.
\end{prop}
\begin{proof} We know that $(q/p)$ is presented by a rank 2 pairing
$V$ of odd type.  Now $V$ contains characteristic elements
$v\in V$.  If one could choose $v$ such that $v\cdot v =\pm 1$, then
$(q/p)$ would be presented by the rank 1 pairing
$v^\perp$, which would be of even type.   This cannot happen in
general, since, in particular, we would need $p$ even.  But in any
case consider $V\oplus \left[
\begin{smallmatrix}0&1\\1&0\end{smallmatrix}\right]$. 
Characteristic elements in
$\left[ \begin{smallmatrix}0&1\\1&0\end{smallmatrix}
\right]$ are of the form $w=2ke_1+2\ell e_2$ and $w\cdot w=4k\ell$,
which can be any multiple of $4$.  Thus in
$V\oplus \left[
\begin{smallmatrix}0&1\\1&0\end{smallmatrix}\right]$ there are
characteristic elements of the form $v+w$ such that
$(v+w)^2=-1$, $0$, $1$, or $2$.  If $(v+w)^2 = \pm 1$, then pass to
the orthogonal complement $(v+w)^\perp$, which is even, of rank 3,
and presents $(q/p)$.  Otherwise, first add on an additional
$\left<+1\right>$ or $\left<-1\right>$ and an additional basis
vector $e$ to $v+w$ to get a  characteristic element
$u= v+w+e$ in $V\oplus \left[
\begin{smallmatrix}0&1\\1&0\end{smallmatrix}\right]\oplus (\pm 1)$,
such that $u^2=\pm 1$.  Passing to $u^\perp$ then provides a rank 4
even pairing presenting $(q/p)$, as required.
\end{proof}

\begin{cor}  If $p$ is a positive integer and $q$ is an integer
prime to $p$, then  any (homology) lens space  $L(p,q)$ admits a topological
embedding in
$\#_n\ S^2\times S^2$, $n\le 4$.
\end{cor}
\begin{proof} We have seen that $L(p,q)$ bounds a simply connected
4-manifold of even type and with $n=b_2\le 4$.  The double of such a
4-manifold is precisely $\#_n\ S^2\times S^2$, by Freedman's
classification of simply connected topological 4-manifolds.
\end{proof}

\section{ Definite Intersection Pairings and Applications}

Here we investigate presenting a linking form by a {\em definite}
intersection pairing.  It is not too hard to do this, but we have to
exert some effort to keep the resulting rank as small as possible. 
We begin with some lemmas useful in applying the blow-up/blow-down
procedure introduced in the preceding section.

\begin{lem}Suppose that $V$ is an odd intersection pairing that is
not positive definite.  Then there is $v\in V$ such that $v\cdot v <
0$ and $v\cdot v$ is odd.
\end{lem}
\begin{proof} There is a $u\in V$ such that $u\cdot u<0$ and there
is $w\in V$ such that $w\cdot w$ is odd.  We only need to consider
further the case that $u\cdot u$ is even and $w\cdot w>0$; otherwise
we are done.  Replacing $w$ by
$-w$ if necessary we may also assume that $u\cdot w\le 0$.

Now let $v=w+ku$, $k\in {\bf Z}$.  We compute that
$$v\cdot v = w\cdot w + 2k u\cdot w + k^2 u\cdot u$$ from which it
is clear that $v\cdot v$ is odd and that for sufficiently large $k$
we also have $v\cdot v < 0$.
\end{proof}

\begin{lem}Suppose that $V$ is an odd intersection pairing that is
not positive definite.  Then there is $u\in V$ such that $u\cdot u =
-n < 0$, $u$ is ordinary (actually 2-divisible), and $n-1$ is not of
the form
$4^a(8b+7)$.
\end{lem}
\begin{proof} There is a $v\in V$ such that $v\cdot v=-(2k+1)<0$ for
some integer $k$.  Set $u=2v$.  Then
$u\cdot u = 4 v\cdot v < 0$.  Also $u$ is ordinary, since there is
some $w\in V$ such that $w\cdot w$ is odd, while
$u\cdot w = 2 v\cdot w$ is even.  Finally, setting
$n=-u\cdot u$, we have $n-1 = 4(2k+1)-1=8k+3$, which is never of the
form
$4^a(8b+7)$.
\end{proof}  

\begin{cor} If $V$ is an intersection pairing of odd type that is
not positive definite, then there is $v\in V\oplus 3\left<+1\right>$
such that $v\cdot v=-1$ and $v$ is ordinary.
\end{cor}
\begin{proof} By the preceding result there is $u\in V$ such that
$u\cdot u = -n<0$, $u$ is ordinary (and, in fact, 2-divisible), and
$n-1$ is not of the form $4^a(8b+7)$.  By number theory (see
\cite{Gr}, for example), $n-1$ can be written as a sum of 3 squares,
$n-1=a_1^2 +a_2^2 + a_3^2$.  Set
$v=u+a_1e_1+a_2e_2+a_3e_3$, where $e_1$, $e_2$, and $e_3$ form a
standard orthonormal basis for $3\left<+1\right>$.  Then
$v\cdot v=u\cdot u + a_1^2 +a_2^2 + a_3^2 = -n+n-1=-1$ and
$v$ is ordinary since there is $w\in V\subset V\oplus
3\left<+1\right>$ such that $w\cdot w$ is odd, but $v\cdot w=u\cdot
w$ is even since $u$ is 2-divisible.
\end{proof}

\begin{prop} The rank 1 linking form $(q/p)$ is presented by a
positive definite intersection pairing of rank $\le 6$ and of odd
type.
\end{prop}
\begin{proof} We know that $(q/p)$ is presented by a rank 2 pairing
$V$ of odd type.  If $V$ happens to be positive definite, then we
are done.  There are two remaining cases, depending on whether $V$
is indefinite or negative definite. 

First suppose that $V$ is indefinite.  Then $V$ contains an ordinary
element $v\in V$ such that $v \cdot v =-n< 0$ where.  If we could
choose $v$ such that $v \cdot v = -1$, then $(q/p)$ would be
presented by the rank 1 form
$v^\perp$.  In general this is impossible to achieve.  We can,
however, blow up the form, passing to $V\oplus 3\left<+1\right>$,
which still represents $(q/p)$.  Here we can find an ordinary $v\in
V\oplus 3\left<+1\right>$ such that $v \cdot v = -1$.  Passing to
$v^\perp$, we obtain a positive definite integral pairing of rank
$2+3-1=4$ presenting
$(q/p)$.

Finally we must consider the case when $(q/p)$ is presented by a
negative definite rank 2 pairing $V$ of odd type. In this situation
we do much the same as before.  We need to blow down twice,
however.  In order to do this, we must blow up $3\left<+1\right>$
{\it twice}.  The net effect is to produce a positive definite
pairing of rank $2+3+3-1-1=6$ and odd type presenting $(q/p)$.  It
should be noted that in any particular case we seem to be able to do
better than this, and we know of no case where the full rank 6
possibility is actually required.
\end{proof}

\begin{rmk} If $(q/p)$ is represented by a rank 2 positive definite
odd pairing, then $(-q/p)$ is represented by a rank 2 negative
definite odd pairing, hence by a rank $\le 6$ positive definite odd
pairing.  (And, of course, if
$(q/p)$ is represented by a rank 2 negative definite odd pairing,
then $(-q/p)$ is represented by a rank 2 positive definite odd
pairing, while $(-q/p)$ is then represented by a rank $\le 6$
positive definite odd pairing.) If
$(q/p)$ is represented by a rank 2 indefinite odd pairing, then
$(-q/p)$ is also represented by a rank 2 indefinite odd pairing,
hence both are represented by rank
$\le 5$ positive definite odd pairings, according to the proof of
the theorem.
\end{rmk}

\begin{cor}  If $p$ is a positive integer and $q$ is an integer
prime to $p$, then  any (homology) lens space  $L(p,q)$ admits a topological
embedding in
$\#_n{\mathbf C}P^2$, for some $n\le 8$.
\end{cor}
\begin{proof}  Replacing $q$ by $-q$ if necessary, we know that
either $L(p,q)$ bounds a simply connected 4-manifold of odd type and
with positive definite intersection pairing and $b_2 = 2$, while
$-L(p,q)$ bounds a simply connected 4-manifold with positive
definite intersection pairing and
$b_2 \le 6$; or $L(p,q)$ and $-L(p,q)$ both bound simply connected
4-manifolds with positive definite intersection pairings of odd type
and $b_2 \le 4$.  In either of the two cases,  the union of the two
4-manifolds along the lens space yields a simply connected
4-manifold, with $b_2\le 8$.  Thus in either case we obtain an
embedding of $L(p,q)$ in a closed, simply connected 4-manifold with
positive definite intersection pairing of odd type and rank at most
8.  It follows from the classification of unimodular intersection
pairings of low rank that the intersection pairing is
diagonalizable.  By Freedman's classification theorem, it only
remains to be sure the Kirby-Siebenmann stable triangulation
obstruction vanishes.  This is taken care of by Remark~\ref{ks}
above, since we can, if necessary, change the Kirby-Siebenmann
invariant of one of the two pieces to make the global invariant
vanish.  Then Freedman's classification theorem shows that we have
$\#_n{\mathbf C}P^2$.
\end{proof}

%
\bibliographystyle{amsplain}

\noindent{Department of Mathematics, Indiana University, Bloomington, IN
47405}

{e-mail: edmonds@indiana.edu}

\end{document}